\documentclass{ifacconf}

\usepackage[utf8]{inputenc}
\usepackage{amsmath}
\usepackage{amssymb}
\usepackage{graphicx}      % include this line if your document contains figures
\usepackage{natbib}        % required for bibliography
\usepackage{color}
\usepackage{array}
\usepackage{bm}
\usepackage[final]{microtype}

\newcommand{\abs}[1]{\left\lvert#1\right\rvert}
\newcommand{\sign}[1]{\text{sign}\left(#1\right)}
\begin{document}
\begin{frontmatter}

\title{Consistent Discretization of a Class of Predefined-Time Stable Systems}

\thanks[footnoteinfo]{This work has been submitted to IFAC for possible publication}

\author[First]{Esteban Jiménez-Rodríguez} 
\author[Second]{Rodrigo Aldana-López}
\author[Third]{Juan D. S\'anchez-Torres}
\author[Fourth,Fifth]{David G\'omez-Guti\'errez}
\author[First]{Alexander G. Loukianov}

\address[First]{Department of Electrical Engineering, Cinvestav-Guadalajara, Jalisco, 45019 M\'exico (e-mail: \{ejimenezr, louk\}@gdl.cinvestav.mx)}
\address[Second]{Department of Computer Science and Systems Engineering, University of Zaragoza, Zaragoza, 50009 España, (e-mail: rodrigo.aldana.lopez@gmail.com)}
\address[Third]{Research Laboratory on Optimal Design, Devices and Advanced Materials -OPTIMA-, Department of Mathematics and Physics, ITESO, Jalisco, 45604 M\'exico (e-mail: dsanchez@iteso.mx)}
\address[Fourth]{Multi-agent autonomous systems lab, Intel Labs, Intel Tecnolog\'ia de M\'exico, Jalisco, 45019 Mexico (e-mail: david.gomez.g@ieee.org)}
\address[Fifth]{Tecnologico de Monterrey, Escuela de Ingenier\'ia y Ciencias, Jalisco, 45138 Mexico}

%Predefined-time stability refers to the property of certain class of fixed-time stable systems whose solutions converge to the origin before an arbitrary user-prescribed time.
\begin{abstract} As the main contribution, this document provides a consistent discretization of a class of fixed-time stable systems, namely predefined-time stable systems. In the unperturbed case, the proposed approach allows obtaining not only a consistent but exact discretization of the considered class of predefined-time stable systems, whereas in the perturbed case, the consistent discretization preserves the predefined-time stability property. All the results are validated through simulations and compared with the conventional explicit Euler scheme, highlighting the advantages of this proposal.
\end{abstract}


\end{frontmatter}
%===============================================================================

\section{Introduction}\label{sec:intro}

Given a system with tunable parameters whose origin is fixed-time stable, generally, it is not straightforward, and sometimes it is even impossible to achieve any desired upper bound of the settling time through the selection of the parameters of the system~\citep[Example 1]{Jimenez-Rodriguez2019}. To overcome this drawback, a class of dynamical systems that exhibit the property of \textit{predefined-time stability} has been studied within the last six years~\citep{Sanchez-Torres2018}.

Predefined-time stability refers to the property that exhibits a particular class of fixed-time stable systems whose solutions converge to the origin before an arbitrary user-prescribed time, which is assigned through an appropriate selection of the tunable parameters of the system~\citep{Sanchez-Torres2018}. Due to its remarkable features, the study of the properties and applications of the predefined-time stability notion has attracted much attention, mainly in continuous time~\citep{Aldana-Lopez2019a,Aldana-Lopez2019,Munoz-Vazquez2019,Sanchez-Torres2019}.

The numerical simulation examples conducted in the above works were done applying the conventional explicit (forward) Euler discretization, with a tiny step size ($1\times10^{-4}$ or $1\times10^{-5}$ time units). However, the right side of the ordinary differential equations exhibiting predefined-time stability does not satisfy the Lipschitz condition at the origin. In this case, the explicit Euler discretization scheme does not guarantee that predefined-time stability property will be preserved~\citep{Levant2013,Huber2016,Polyakov2019}. In other words, the discrete-time equation solutions may be inconsistent with the solutions of the continuous-time one.

In this sense, this work concerns with the discretization of a class of predefined-time stable systems. The discretization process uses the idea of~\cite[Figure 1.1]{Polyakov2019}, but different from the contribution of the mentioned paper, which achieves consistent discretizations of \textit{homogeneous finite-time and practically fixed-time stable systems}, the proposed approach allows the consistent discretization of the class of predefined-time stable systems proposed by~\cite{Aldana-Lopez2019,Jimenez-Rodriguez2019}. Hence, this document provides a consistent discretization of a class of fixed-time stable systems.

The rest of this paper is organized as follows: $\mathcal{K}^1$ functions and the considered class of predefined-time stable systems are introduced in Section~\ref{sec:prelim}. Then, the exact discretization of the considered class of predefined-time stable systems is developed in Section~\ref{sec:exact}, in two different ways. Finally, Section~\ref{sec:simu} considers a consistent discretization of perturbed systems with discontinuous predefined-time control.

%First, discrete-time equivalent systems of a class of continuous-time predefined-time stable systems are obtained; we say that they are equivalent in the sense that their solutions are equal to the continuous-time system solutions at the sampling instants. Then, a consistent discretization of predefined-time controlled perturbed systems is introduced; we say that it is consistent since it preserves the predefined-time stability property of the closed-loop system. All the results are illustrated with numerical examples, and are compared with explicit Euler discretization schemes to show the key differences and advantages.

\section{Preliminaries}\label{sec:prelim}

\subsection{Class~$\mathcal{K}^1$ functions}

The following definition was inspired on the comparison class $\mathcal{K}$ functions~\citep[Definition 1]{Kellett2014}

\begin{defn}[$\mathcal{K}^1$ functions] \label{def:class_K1} A scalar continuous function $\kappa: \mathbb{R}_{\geq 0} \to \left[ 0, 1 \right)$ is said to belong to class $\mathcal{K}^1$, denoted as $\kappa \in \mathcal{K}^1$, if it is strictly increasing, $\kappa(0) = 0$ and $\lim_{r \to \infty} \kappa(r) = 1$.
\end{defn}

\begin{rem}[Invertibility of class $\mathcal{K}^1$ functions] Let $\kappa \in \mathcal{K}^1$.
\begin{itemize}
\item[\textit{(i)}] $\kappa$ is injective since it is strictly increasing;
\item[\textit{(ii)}] $\kappa$ is onto since it is continuous, $\kappa(0) = 0$ and $\lim_{r \to \infty} \kappa(r) = 1$.
\end{itemize}
Hence, $\kappa$ is bijective and its inverse exist.
\end{rem}

\begin{rem}[A connection with probability theory] Let $\kappa \in \mathcal{K}^1$ be continuously differentiable on its domain. In this case, there exists a continuous function $\Phi:\mathbb{R}_{\geq 0}\to\mathbb{R}_+$ such that $\kappa(r)=\int_{0}^{r}\Phi(z)\text{d}z$. Moreover, since $\lim_{r \to \infty} \kappa(r) = 1$, the function $\Phi$ is required to satisfy $\int_0^\infty \Phi(z)\text{d}z=1$, i.e., functions $\Phi$ and $\kappa$ can be viewed as probability density functions and cumulative distribution functions, respectively, of positive random variables.
\end{rem}

\begin{exmp} Examples of class $\mathcal{K}^1$ functions are:
\begin{itemize}
\item $\kappa(r)=\frac{2}{\pi}\arctan(ar)$, with $a>0$;
\item $\kappa(r)=\frac{r}{r+a}$, with $a>0$;
\item $\kappa(r)=1-a^{-r}$, with $a>1$;
\item $\kappa(r)=P(a,r)$, with $a>0$, where  $P(a,r)=\frac{\gamma(a,r)}{\Gamma(a)}$ is the regularized Incomplete Gamma Function, $\gamma(a,r)=\int_0^r t^{a-1}e^{-t}\text{d}t$ is the Incomplete Gamma Function and $\Gamma(a)=\int_0^\infty t^{a-1}e^{-t}\text{d}t$ is the Gamma Function~\citep{Abramowitz1964};
\item $\kappa(r)=I(a_1,a_2,\frac{r}{r+1})$, with $a_1,a_2>0$, where $I(a_1,a_2,r)=\frac{b(a_1,a_2,r)}{\mathcal{B}(a_1,a_2)}$ is the regularized Incomplete Beta Function, $b(a_1,a_2,r)=\int_0^{r} t^{a_1-1}(1-t)^{a_2-1}\text{d}t$ is the Incomplete Beta Function and $\mathcal{B}(a_1,a_2)=\int_0^1 t^{a_1-1}(1-t)^{a_2-1}\text{d}t$ is the Beta Function~\citep{Abramowitz1964}.
\end{itemize}
\end{exmp}

\subsection{A class of predefined-time stable systems}\label{subsec:ptsys}

Consider (the scalar form of) the class of systems introduced by~\cite{Aldana-Lopez2019} and~\cite{Jimenez-Rodriguez2019}
\begin{equation} \label{eq:ptsys}
\dot{x} = - \frac{1}{\rho_1 (1-\rho_2)} \frac{\kappa(\abs{x})^{\rho_2}}{\kappa'(\abs{x})}\sign{x}, \ \ x(0)=x^0,
\end{equation}
where $x:\mathbb{R}_{\geq 0}\to\mathbb{R}$ is the state of the system, $\rho_1 > 0$ and $0 \leq \rho_2 < 1$ are tunable parameters, $\kappa \in \mathcal{K}^1$ is continuously differentiable, and $\sign{\cdot}$ stands for the signum function given by
\[
\sign{x}\in\left\lbrace
\begin{array}{ccc}
    \{1\} & \text{if} & x>0 \\
    \{-1\} & \text{if} & x<0 \\
    {[-1,1]} & \text{if} & x=0.
\end{array}\right.
\]

From the theory of ordinary differential equations, \eqref{eq:ptsys} is a separable first-order equation, whose solution $x(t)$ satisfies
\[
\int_{x^0}^{x(t)}(1-\rho_2)\frac{\kappa'(\abs{\xi})}{\kappa(\abs{\xi})^{\rho_2}}\sign{\xi}\mathrm{d}\xi=-\int_{0}^{t}\frac{1}{\rho_1}\mathrm{d}\tau.
\]
Applying the change of variable $\xi \to \kappa(\abs{\xi})^{1-\rho_2}\sign{\xi}$ in the left side integral and integrating both sides of the above equation, it yields
\[
\kappa(\abs{x(t)})^{1-\rho_2}=\kappa(\abs{x^0})^{1-\rho_2}-\frac{t}{\rho_1},
\]
for $\kappa(\abs{x^0})^{1-\rho_2}-\frac{t}{\rho_1} \geq 0$. From the above, the solution of~\eqref{eq:ptsys} is
\begin{multline}\label{eq:ptsol}
x(t)=\\
\small{
\left\lbrace
\begin{array}{l}
\kappa^{-1}\left(\left[\kappa\left(\abs{x^0}\right)^{1-\rho_2}-\frac{t}{\rho_1}\right]^{\frac{1}{1-\rho_2}}\right)\sign{x^0} \\
\qquad\qquad\qquad\qquad\qquad\qquad\qquad \text{ if } 0\leq t \leq \rho_1\kappa\left(\abs{x^0}\right)^{1-\rho_2} \\
0 \\
\qquad\qquad\qquad\qquad\qquad\qquad\qquad \text{ if } t>\rho_1\kappa\left(\abs{x^0}\right)^{1-\rho_2},
\end{array}
\right.
}
\end{multline}
where $\kappa^{-1}$ stands for the inverse of the function $\kappa\in\mathcal{K}$.
% \begin{equation}\label{eq:ptsol}
% \phi(t,x^0)=\kappa^{-1}\left(\left[\kappa(\abs{x^0})^{1-\rho_2}-\frac{t}{\rho_1}\right]^{\frac{1}{1-\rho_2}}\right)\sign{x^0}
% \end{equation}
% for $0\leq t \leq \rho_1\kappa\left(\abs{x^0}\right)^{1-\rho_2}$, and $\phi(t,x^0)=0$ for $t>\rho_1\kappa\left(\abs{x^0}\right)^{1-\rho_2}$.

Hence, one can conclude that the origin $x=0$ of~\eqref{eq:ptsys} is finite-time stable and the settling-time function is $T(x^0)=\rho_1\kappa\left(\abs{x^0}\right)^{1-\rho_2}$. But the settling-time function $T(x^0)=\rho_1\kappa\left(\abs{x^0}\right)^{1-\rho_2}<\rho_1$ is bounded, then the origin $x=0$ of~\eqref{eq:ptsys} is, in fact, fixed-time stable. 

Moreover, for any desired upper bound of the settling-time function $T_c$, there exists an appropriate selection of the parameter $\rho_1=T_c$ such that the convergence time will always be less than $T_c$ for any initial condition $x^0$. This property is known as \textbf{predefined-time stability}~(see Appendix~\ref{ap:pts}).

\begin{rem} It is worth noticing that the right-hand side of~\eqref{eq:ptsys} is continuous and non-Lipschitz if $0 < \rho_2 < 1$, and discontinuous if $\rho_2=0$. In any case, all the above derivation remains true.
\end{rem}

\section{Exact discretization of the considered class of predefined-time stable systems} \label{sec:exact}

In this section, we develop a discrete-time system whose solutions are equal to the solution of~\eqref{eq:ptsys} in the sampling instants. Hence, we provide the exact form of simulating the class of predefined-time stable systems represented by~\eqref{eq:ptsys}.

\subsection{Direct construction}

Consider the continuous-time system~\eqref{eq:ptsys} and let $x_k = x(kh)$, where $h>0$ is the sampling step size, and $k\in\mathbb{N}$.

To obtain the exact discretization of system~\eqref{eq:ptsys}, one can apply the same procedure as in Subsection~\ref{subsec:ptsys}, integrating the left side between the samples $x_{k}$ and $x_{k+1}$, and the right side between $kh$ and $(k+1)h$, respectively, obtaining
\begin{equation}\label{eq:exact}
\int_{x_{k}}^{x_{k+1}}(1-\rho_2)\frac{\kappa'(\abs{\xi})}{\kappa(\abs{\xi})^{\rho_2}}\sign{\xi}\mathrm{d}\xi=-\int_{kh}^{(k+1)h}\frac{1}{\rho_1}\mathrm{d}\tau.
\end{equation}

Recall the transformation
\begin{equation}\label{eq:transformation}
y_{\rho_2}  \colon  \begin{array}{>{\displaystyle}r @{} >{{}}c<{{}} @{} >{\displaystyle}l} 
                        \mathbb{R} &\rightarrow& (-1,1) \\ 
                        x &\mapsto& y_{\rho_2}(x) = \kappa(\abs{x})^{1-\rho_2}\sign{x}
                    \end{array}
\end{equation}
used in Subsection~\ref{subsec:ptsys} to find the solutions of~\eqref{eq:ptsys}. This transformation is invertible, and its inverse $y_{\rho_2}^{-1}:(-1, 1)\to\mathbb{R}$ is given by $$y_{\rho_2}^{-1}(x)=\kappa^{-1}\left(\abs{x}^{\frac{1}{1-\rho_2}}\right)\sign{x}.$$

Using~\eqref{eq:transformation} in the left side integral of~\eqref{eq:exact}, one obtains
\[
\int_{\kappa(\abs{x_{k}})^{1-\rho_2}\sign{x_{k}}}^{\kappa(\abs{x_{k+1}})^{1-\rho_2}\sign{x_{k+1}}}\sign{y}\mathrm{d}y=-\int_{kh}^{(k+1)h}\frac{1}{\rho_1}\mathrm{d}\tau.
\]

Hence, integrating both sides of the above expression and having in mind that $\sign{x_{k+1}}=\sign{x_{k}}$ for the first-order system~\eqref{eq:ptsys}, the above equation turns into
\begin{multline}\label{eq:implicit}
\kappa(\abs{x_{k+1}})^{1-\rho_2}\sign{x_{k+1}} \\= \kappa(\abs{x_{k}})^{1-\rho_2}\sign{x_{k}} - \frac{h}{\rho_1}\sign{x_{k+1}}.
\end{multline}

Finally, let $\alpha\in\mathbb{R}_+$ and consider the discontinuous function $\theta_{\rho_2,\alpha}:\mathbb{R}\to(-\alpha-1,\alpha+1)$ defined by $\theta_{\rho_2,\alpha}(x)=\left(\kappa(\abs{x})^{1-\rho_2}+\alpha\right)\sign{x}$. The inverse of $\theta_{\rho_2,\alpha}$ is the continuous function $\theta_{\rho_2,\alpha}^{-1}:(-\alpha-1,\alpha+1)\to\mathbb{R}$ given by
\begin{equation}\label{eq:inverse}
\theta_{\rho_2,\alpha}^{-1}(x)=\kappa^{-1}\left(\left(\max\{\abs{x}-\alpha,0\}\right)^{\frac{1}{1-\rho_2}}\right)\sign{x}
\end{equation}

Thus, replacing~\eqref{eq:inverse} in~\eqref{eq:implicit}, it yields
\begin{align}\label{eq:discrete}
\begin{split}
x_{k+1} &=\theta_{\rho_2,h/\rho_1}^{-1}\left(\kappa(\abs{x_{k}})^{1-\rho_2}\sign{x_{k}}\right),\\
        &=\kappa^{-1}\left(\max\left\{\kappa(\abs{x_{k}})^{1-\rho_2}-\frac{h}{\rho_1},0\right\}^{\frac{1}{1-\rho_2}}\right)\sign{x_k},
\end{split}
\end{align}
with $x_0=x^0$.

This following theorem summarizes the previous constructive result:

\begin{thm}\label{thm:equivalence} The solution of the discrete-time system~\eqref{eq:discrete} is equivalent to the solution~\eqref{eq:ptsol} of the continuous-time system~\eqref{eq:ptsys} at the sample instants. This is, $x_k = x(kh)$, for every $k\in\mathbb{N}$.
\end{thm}

Theorem~\ref{thm:equivalence} has the following attractive corollary:

\begin{cor} For any $x^0\in\mathbb{R}$, the solution of the discrete-time system~\eqref{eq:discrete} satisfies $x_k=0, \text{ for } k\geq \left\lceil \frac{\rho_1}{h} \right\rceil,$ where $\left\lceil\cdot\right\rceil$ stands for the ceiling function.
\end{cor}
\begin{pf} The solution of the continuous-time system~\eqref{eq:ptsys} satisfies
\[
x(t)=0, \text{ for } t\geq \rho_1\kappa\left(\abs{x^0}\right)^{1-\rho_2}.
\]
Moreover, from Theorem~\ref{thm:equivalence}, we can establish for the solution of the discrete-time system~\eqref{eq:discrete} that
\[
x_k=x(kh)=0, \text{ for } kh\geq \rho_1\kappa\left(\abs{x^0}\right)^{1-\rho_2}.
\]
which also holds for $kh\geq \rho_1$, given that $\kappa\in\mathcal{K}$.

Hence, the result follows.\hfill$\blacksquare$
\end{pf}

%\textcolor{red}{Simulation example here.}

\begin{exmp} Consider system~\eqref{eq:ptsys} with the particular selections of $\kappa(r)=\frac{2}{\pi}\arctan(r)$, $\rho_1=1$, $\rho_2=0.5$ and $x^0=10$. 

Figures~\ref{fig:fig1} and~\ref{fig:fig2} present the comparison of the exact discretization~\eqref{eq:discrete} with the conventional explicit Euler discretization for two different time steps. In both cases, the solution of the discrete-time system~\eqref{eq:discrete} accurately recovers the exact solution~\eqref{eq:ptsol} at the sampling instants. For $h=0.02$ (Fig.~\ref{fig:fig1}), one can note small differences between the exact solution and the explicit Euler discretization, whereas, for $h=0.0651318636497$ (Fig.~\ref{fig:fig2}), the explicit Euler discretization becomes unstable.

\begin{figure}[ht]
\begin{center}
\includegraphics[width=8.4cm]{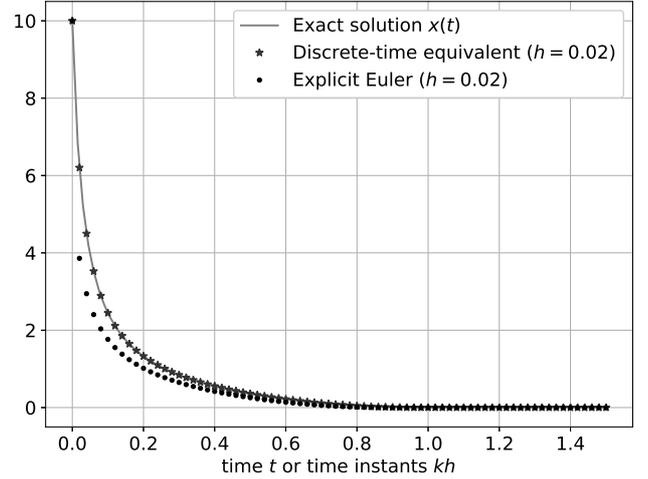}
\caption{Comparison of the exact solution (solid line), the discrete-time equivalent (star points) and the explicit Euler discretization (round points) with a step size of $h=0.02$.}
\label{fig:fig1}
\end{center}
\end{figure}

\begin{figure}[ht]
\begin{center}
\includegraphics[width=8.4cm]{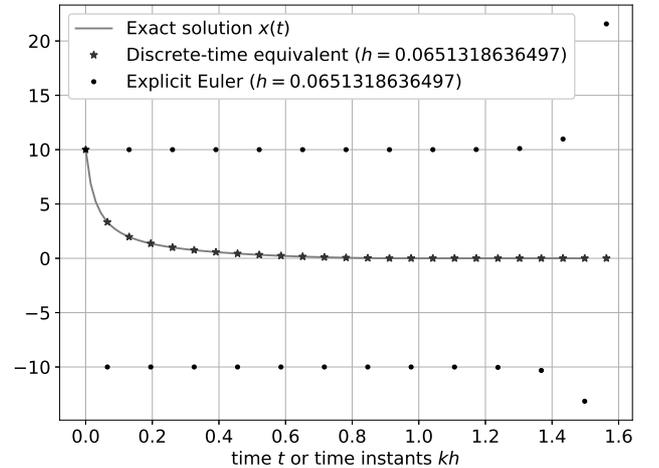}
\caption{Comparison of the exact solution (solid line), the discrete-time equivalent (star points) and the explicit Euler discretization (round points) with a step size of $h=0.0651318636497$.}
\label{fig:fig2}
\end{center}
\end{figure}
\end{exmp}

\subsection{Parallel with implicit Euler discretization}

Consider the nonlinear transformation given by~\eqref{eq:transformation}, and define $w(t)=y_{\rho_2}(x(t))$. Its time-derivative along the trajectories of system~\eqref{eq:ptsys} is
\begin{equation}\label{eq:trans_sys}
    \dot{w} = -\frac{1}{\rho_1}\sign{w}, \ \ w(0)=y_{\rho_2}(x^0),
\end{equation}
since $\sign{w}=\sign{x}$.

\begin{rem} If the \textit{explicit or forward Euler procedure} with time step $h>0$ were used to simulate system~\eqref{eq:trans_sys}, it would lead to oscillations with amplitude $\mathcal{O}(h)$ in the neighborhood of the origin $w=0$~\citep{Drakunov1989,Utkin1994}, and these oscillations would also be translated into oscillations in the original variable $x$. Those mentioned oscillations, known as numerical chattering, are an undesired effect of an inadequate discretization, since the solution~\eqref{eq:ptsol} does not oscillate.
\end{rem}

Using the \textit{implicit or backward Euler method} (as in~\citep{Drakunov1989,Utkin1994}) with time step $h>0$ for the discretization of the transformed system~\eqref{eq:trans_sys}, it yields
\[
w_{k+1} = w_{k} - \frac{h}{\rho_1}\sign{w_{k+1}},
\]
where $w_{k}=y_{\rho_2}(x_{k})$ and $x_{k}=x(kh)$. Rewriting the above in terms of $x_k$, one obtains~\eqref{eq:exact}, as before.

\section{Consistent discretization of predefined-time control of perturbed systems}\label{sec:simu}

\subsection{Continuous-time predefined-time stabilization of first order perturbed systems}
%\subsection{Continuous-time predefined-time stabilization of perturbed systems}

Consider the perturbed control system
\begin{equation}\label{eq:control_sys}
    \dot{x} = u + \Delta(t,x),\,\, x(0)=x^0,
\end{equation}
where $x:\mathbb{R}_{\geq 0}\to\mathbb{R}$ is the state of the system, $u\in\mathbb{R}$ is the control input signal, and $\Delta:\mathbb{R}_{\geq 0}\times\mathbb{R}\to\mathbb{R}$ is an unknown perturbation term which is assumed to be bounded of the form $\sup_{(t,x)\in\mathbb{R}_{\geq 0}\times\mathbb{R}} \abs{\Delta(t,x)} \leq \delta$, with $\delta\in\mathbb{R}_+$ known.

It is well known that, in continuous time, the perturbation term $\Delta(t,x)$ can only be entirely rejected by a discontinuous control term, since no conditions of smoothness, Lipschitz continuity neither continuity are assumed~\citep{Jimenez-Rodriguez2019}.

According to the above, a suitable feedback predefined-time controller design is
\begin{equation}\label{eq:control_in}
    u = - \left(\frac{1}{\rho_1} + \rho_3\kappa'(0)\right) \frac{1}{\kappa'(\abs{x})}\sign{x},
\end{equation}
where $\rho_1>0$, $\rho_3\geq\delta$, and $\kappa\in\mathcal{K}^{1}$ is continuously differentiable and such that $\kappa'(0):=\kappa'(0^{+}) \geq \kappa'(r)$ for all $r\in\mathbb{R}_{\geq 0}$. Note that~\eqref{eq:control_in} is an augmented gain version of the right side of~\eqref{eq:ptsys} with $\rho_2=0$.

To verify the predefined-time convergence of the closed-loop system~\eqref{eq:control_sys}-\eqref{eq:control_in}, note that
\begin{align*}
    \frac{d\abs{x}}{dt} & = - \left(\frac{1}{\rho_1} + \rho_3\kappa'(0)\right) \frac{1}{\kappa'(\abs{x})} + \Delta(t, x)\sign{x}\\
    & \leq - \frac{1}{\rho_1\kappa'(\abs{x})} - \rho_3\frac{\kappa'(0)}{\kappa'(\abs{x})} + \abs{\Delta(t, x)}\\
    & \leq - \frac{1}{\rho_1\kappa'(\abs{x})} - (\rho_3 - \delta)\\
    & \leq - \frac{1}{\rho_1\kappa'(\abs{x})},
\end{align*}
for $x\neq 0$. Reasoning as in Subsection~\ref{subsec:ptsys}, and using the comparison lemma in the above differential inequality, one can easily conclude that $\abs{x(t)}=x(t)=0$ for $t\geq \rho_1$.

\subsection{Consistent discretization/simulation}

Consider transformation~\eqref{eq:transformation}, with $\rho_2=0$, given by $z(t)=y_{0}(x(t))=\kappa(\abs{x(t)})\sign{x(t)}$. Its derivative along the trajectories of the closed-loop system~\eqref{eq:control_sys}-\eqref{eq:control_in} is
\begin{equation}\label{eq:trans_control_sys}
    \dot{z} = -\beta\sign{z}+f(t,x),\,\, z(0)=y_0(x^0),
\end{equation}
where $\beta:=\frac{1}{\rho_1}+\rho_3\kappa'(0)$, and $f(t,x)=\kappa'(\abs{x})\Delta(t,x)$ is the perturbation term in the transformed coordinate $z$, which complies to $\sup_{(t,x)\in\mathbb{R}_{\geq 0}\times\mathbb{R}} \abs{f(t,x)} \leq \kappa'(0)\delta$.

Again, using the \textit{implicit Euler method} with time step $h>0$ for the discretization of the system~\eqref{eq:trans_control_sys}, it yields
\begin{equation}\label{eq:implicit_cont}
    z_{k+1} = z_{k} - h\beta\sign{z_{k+1}} + hf_k,
\end{equation}

where $f_k=f(kh,x_k)$ is the perturbation term at the sample instants.

Hence, replacing the function $\theta_{0,h\beta}^{-1}$~\eqref{eq:inverse} into~\eqref{eq:implicit_cont}, one obtains
\begin{multline}\label{eq:discrete_cont}
    x_{k+1}=\\\kappa^{-1}\left(\max\left\{\abs{\kappa(\abs{x_{k}})\sign{x_k}+hf_k}-h\beta,0\right\}\right)\times\\\sign{\kappa(\abs{x_{k}})\sign{x_k}+hf_k},
\end{multline}
with $x_0=x^0$.

The discrete-time system~\eqref{eq:discrete_cont} provides a way of simulating the closed-loop continuous-time system~\eqref{eq:control_sys}-\eqref{eq:control_in}. Moreover, it preserves the predefined-time stability behavior, as stated in the following proposition:

\begin{prop} For any $x^0\in\mathbb{R}$, the solution of the discrete-time system~\eqref{eq:discrete_cont} satisfies $x_k=0, \text{ for } k\geq \left\lceil \frac{\rho_1}{h} \right\rceil$.
\end{prop}
\begin{pf} Taking absolute value in both sides of~\eqref{eq:discrete_cont} and using the triangle inequality one gets,
\begin{align*}
\abs{x_{k+1}}   & = \kappa^{-1}\left(\max\left\{\abs{\kappa(\abs{x_{k}})\sign{x_k}+hf_k}-h\beta,0\right\}\right)\\
                & \leq \kappa^{-1}\left(\max\left\{\kappa(\abs{x_{k}})+h\abs{f_k}-h\beta,0\right\}\right).
\end{align*}
On the other hand, from $\beta:=\frac{1}{\rho_1}+\rho_3\kappa'(0)$, $\rho_3\geq\delta$, $\abs{f_k}\leq \kappa'(0)\delta$ and the above, the inequality
\[
\abs{x_{k+1}}\leq \kappa^{-1}\left(\max\left\{\kappa(\abs{x_{k}})-\frac{h}{\rho_1},0\right\}\right)
\]
is obtained.

Finally, noticing that the right side of the above inequality is equal to the absolute value of the right side of~\eqref{eq:discrete} with $\rho_2=0$, and that it is a non-decreasing function of $x_k$, we get
\[
\abs{x_k} \leq 0,\,\, \forall k\geq \left\lceil \frac{\rho_1}{h} \right\rceil
\]
applying Lemma~\ref{lem:comparison}.\hfill$\blacksquare$
\end{pf}

\begin{exmp} Consider the closed-loop system~\eqref{eq:control_sys}-\eqref{eq:control_in} with the particular selections of $\kappa(r)=\frac{2}{\pi}\arctan(r)$, $\rho_1=1$, $\rho_3=1.1$ and $x^0=10$. It is assumed that $\Delta(t,x)=\sin{10\pi t}$.

Figure~\ref{fig:fig3} present the comparison of the consistent discretization~\eqref{eq:discrete_cont} with the conventional explicit Euler discretization for $h=0.05$. One can see that $x_k=0$ for $k\geq\left\lceil\frac{\rho_1}{h}\right\rceil$ for the consistent discretization~\eqref{eq:discrete_cont}, whereas the explicit Euler discretization induce undesired oscillations around the origin.

\begin{figure}[ht]
\begin{center}
\includegraphics[width=8.4cm]{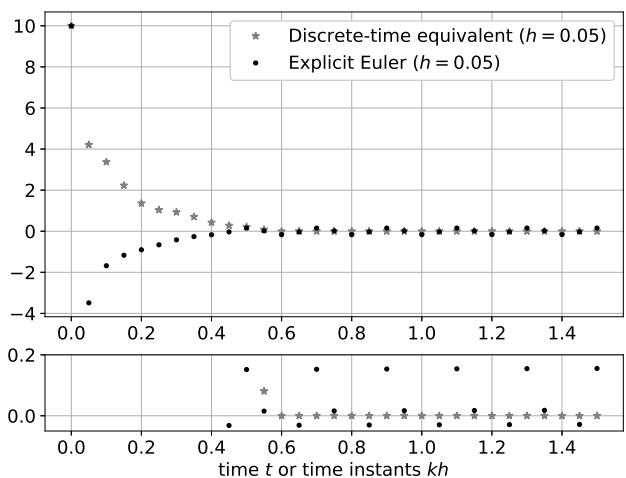}
\caption{Comparison of the the discrete-time equivalent (star points) and the explicit Euler discretization (round points) with a step size of $h=0.05$ for the predefined-time control of perturbed systems.}
\label{fig:fig3}
\end{center}
\end{figure}
\end{exmp}

\section{Conclusion}
This paper presented the development of a consistent discretization for the class of predefined-time stable systems introduced in~\cite{Aldana-Lopez2019,Jimenez-Rodriguez2019}.

The proposed approach allowed the exact discretization of the considered class of systems when no perturbations were assumed. In the perturbed case, the developed consistent discretization preserved the predefined-time stability property.

All the results were confirmed through numerical simulations and compared with the conventional explicit Euler scheme. Even with relatively large time steps, the proposed discretization worked well, as expected, whereas the explicit Euler discretization produced unstable oscillations.

\begin{ack}
Esteban Jim\'enez acknowledges to CONACYT--M\'exico for the D.Sc. scholarship number 481467 and the project 252405.
\end{ack}

%\bibliography{ifacconf}             % bib file to produce the bibliography

\begin{thebibliography}{15}
\providecommand{\natexlab}[1]{#1}
\providecommand{\url}[1]{\texttt{#1}}
\providecommand{\urlprefix}{URL }
\expandafter\ifx\csname urlstyle\endcsname\relax
  \providecommand{\doi}[1]{doi:\discretionary{}{}{}#1}\else
  \providecommand{\doi}{doi:\discretionary{}{}{}\begingroup
  \urlstyle{rm}\Url}\fi

\bibitem[{Abramowitz and Stegun(1965)}]{Abramowitz1964}
Abramowitz, M. and Stegun, I.A. (1965).
\newblock \emph{Handbook of Mathematical Functions: With Formulas, Graphs, and
  Mathematical Tables}.
\newblock Dover Publications, 9 edition.

\bibitem[{{Aldana-L{\'o}pez} et~al.(2019){Aldana-L{\'o}pez},
  {G{\'o}mez-Guti{\'e}rrez}, {Jim{\'e}nez-Rodr{\'{\i}}guez},
  {S{\'a}nchez-Torres}, and {Defoort}}]{Aldana-Lopez2019a}
{Aldana-L{\'o}pez}, R., {G{\'o}mez-Guti{\'e}rrez}, D.,
  {Jim{\'e}nez-Rodr{\'{\i}}guez}, E., {S{\'a}nchez-Torres}, J.D., and
  {Defoort}, M. (2019).
\newblock {Enhancing the settling time estimation of a class of fixed-time
  stable systems}.
\newblock \emph{International Journal of Robust and Nonlinear Control}, 29(12),
  4135--4148.

\bibitem[{Aldana-López et~al.(2019)Aldana-López, Gómez-Gutiérrez,
  Jiménez-Rodríguez, Sánchez-Torres, and Defoort}]{Aldana-Lopez2019}
Aldana-López, R., Gómez-Gutiérrez, D., Jiménez-Rodríguez, E.,
  Sánchez-Torres, J.D., and Defoort, M. (2019).
\newblock On the design of new classes of fixed-time stable systems with
  predefined upper bound for the settling time.

\bibitem[{Bitsoris and Gravalou(1995)}]{Bitsoris1995}
Bitsoris, G. and Gravalou, E. (1995).
\newblock Comparison principle, positive invariance and constrained regulation
  of nonlinear systems.
\newblock \emph{Automatica}, 31(2), 217 -- 222.
\newblock \doi{https://doi.org/10.1016/0005-1098(94)E0044-I}.

\bibitem[{Drakunov and Utkin(1989)}]{Drakunov1989}
Drakunov, S. and Utkin, V. (1989).
\newblock On discrete-time sliding modes.
\newblock \emph{IFAC Proceedings Volumes}, 22(3), 273 -- 278.
\newblock \doi{https://doi.org/10.1016/S1474-6670(17)53647-2}.

\bibitem[{Filippov(1988)}]{Filippov1988}
Filippov, A.F. (1988).
\newblock \emph{Differential equations with discontinuous righthand sides}.
\newblock Kluwer Academic Publishers Group, Dordrecht.

\bibitem[{{Huber} et~al.(2016){Huber}, {Acary}, and {Brogliato}}]{Huber2016}
{Huber}, O., {Acary}, V., and {Brogliato}, B. (2016).
\newblock Lyapunov stability and performance analysis of the implicit discrete
  sliding mode control.
\newblock \emph{IEEE Transactions on Automatic Control}, 61(10), 3016--3030.
\newblock \doi{10.1109/TAC.2015.2506991}.

\bibitem[{Jim{\'e}nez-Rodr{\'i}guez et~al.(2019)Jim{\'e}nez-Rodr{\'i}guez,
  Mu\~noz V\'azquez, S{\'a}nchez-Torres, Defoort, and
  Loukianov}]{Jimenez-Rodriguez2019}
Jim{\'e}nez-Rodr{\'i}guez, E., Mu\~noz V\'azquez, A.J., S{\'a}nchez-Torres,
  J.D., Defoort, M., and Loukianov, A.G. (2019).
\newblock {A Lyapunov-like Characterization of Predefined-Time Stability}.
\newblock \emph{ArXiv e-prints}.
\newblock \urlprefix\url{https://arxiv.org/abs/1910.14604}.

\bibitem[{Kellett(2014)}]{Kellett2014}
Kellett, C.M. (2014).
\newblock A compendium of comparison function results.
\newblock \emph{Mathematics of Control, Signals, and Systems}, 26(3), 339--374.
\newblock \doi{10.1007/s00498-014-0128-8}.

\bibitem[{{Levant}(2013)}]{Levant2013}
{Levant}, A. (2013).
\newblock On fixed and finite time stability in sliding mode control.
\newblock In \emph{52nd IEEE Conference on Decision and Control}, 4260--4265.

\bibitem[{{Mu\~noz-V\'azquez} et~al.(2019){Mu\~noz-V\'azquez},
  {S\'anchez-Torres}, {Jim\'enez-Rodr\'iguez}, and
  {Loukianov}}]{Munoz-Vazquez2019}
{Mu\~noz-V\'azquez}, A.J., {S\'anchez-Torres}, J.D., {Jim\'enez-Rodr\'iguez},
  E., and {Loukianov}, A. (2019).
\newblock Predefined-time robust stabilization of robotic manipulators.
\newblock \emph{IEEE/ASME Transactions on Mechatronics}, 1--1.
\newblock \doi{10.1109/TMECH.2019.2906289}.

\bibitem[{Polyakov et~al.(2019)Polyakov, Efimov, and Brogliato}]{Polyakov2019}
Polyakov, A., Efimov, D., and Brogliato, B. (2019).
\newblock Consistent discretization of finite-time and fixed-time stable
  systems.
\newblock \emph{SIAM Journal on Control and Optimization}, 57(1), 78--103.
\newblock \doi{10.1137/18M1197345}.

\bibitem[{{S\'anchez-Torres} et~al.(2019){S\'anchez-Torres}, {Defoort}, and
  {Mu\~noz-V\'azquez}}]{Sanchez-Torres2019}
{S\'anchez-Torres}, J.D., {Defoort}, M., and {Mu\~noz-V\'azquez}, A.J. (2019).
\newblock Predefined-time stabilisation of a class of nonholonomic systems.
\newblock \emph{International Journal of Control}, 1--8.
\newblock \doi{10.1080/00207179.2019.1569262}.

\bibitem[{S\'anchez-Torres et~al.(2018)S\'anchez-Torres, G\'omez-Guti\'errez,
  L\'opez, and Loukianov}]{Sanchez-Torres2018}
S\'anchez-Torres, J.D., G\'omez-Guti\'errez, D., L\'opez, E., and Loukianov,
  A.G. (2018).
\newblock A class of predefined-time stable dynamical systems.
\newblock \emph{IMA Journal of Mathematical Control and Information}, 35(Suppl
  1), i1--i29.
\newblock \doi{10.1093/imamci/dnx004}.

\bibitem[{Utkin(1994)}]{Utkin1994}
Utkin, V.I. (1994).
\newblock Sliding mode control in discrete-time and difference systems.
\newblock In A.~Zinober (ed.), \emph{Variable Structure and Lyapunov Control},
  chapter~5, 87--107. Springer, Berlin Heidelberg.

\end{thebibliography}
                                    % with bibtex (preferred)

\appendix

\section{Predefined-time stability}\label{ap:pts}

Predefined-time stability refers to the property that exhibits a particular class of fixed-time stable systems with tunable parameters, for which an upper bound of the settling-time function can be arbitrarily chosen through a suitable selection of the parameters~\citep[Defintion~2]{Jimenez-Rodriguez2019}.

This notion is formally defined considering an autonomous system of the form
\begin{equation} \label{eq:sys}
\dot{\bm{x}} = \bm{f}(\bm{x};\bm{\rho}), \ \ \bm{x}(0)=\bm{x}_0,
\end{equation}
with $\bm{x}:\mathbb{R}_{\geq 0}\to\mathbb{R}^n$ the system state, $\bm{\rho}\in\mathbb{R}^l$ the \textit{tunable} parameters of~\eqref{eq:sys}, and $\bm{f}:\mathbb{R}^n\rightarrow\mathbb{R}^n$ a nonlinear function. 

For the origin $\bm{x}=\bm{0}$ to be a finite-, in particular, a predefined-time stable equilibrium of~\eqref{eq:sys}, the function $\bm{f}$ must be a non Lipschitz (maybe discontinuous) function of $\bm{x}$. Then, $\bm{f}$ is assumed to be such that the solutions of~\eqref{eq:sys} exist and are unique in the sense of Filippov~\citep{Filippov1988}.

\begin{defn}\label{def:predefined} The origin of~\eqref{eq:sys} is said to be \textbf{predefined-time stable} if it is fixed-time stable and for any $T_c\in\mathbb{R}_+$, there exists some $\bm{\rho}\in\mathbb{R}^l$ such that the settling-time function of~\eqref{eq:sys} satisfies \[\sup_{\bm{x}_0\in\mathbb{R}^n}T(\bm{x}_0)\leq T_c.\]
\end{defn}

\section{A comparison lemma for discrete-time systems}\label{ap:discrete}
The following is a particular case of \citep[Proposition 1]{Bitsoris1995}.

\begin{lem}\label{lem:comparison} Let $u_k$ be the solution of
\[
u_{k+1}=f(u_k),\,\, u_0=u^0, \,\, u\in\mathbb{R},
\]
where $f:\mathbb{R}\to\mathbb{R}$ is continuous and non-decreasing, and let $v_k$ be such that
\[
v_{k+1}\leq f(v_k),\,\, v_0\leq u^0, \,\, v\in\mathbb{R}.
\]
Then, $v_k\leq u_k$, for all $k\in\mathbb{N}$.
\end{lem}
\begin{pf} We proceed by induction. The base case follows from the hypothesis
\[
v_0\leq u^0 = u_0.
\]
Now, assume that $v_k\leq u_k$. Then, 
\[
v_{k+1} \leq f(v_k) \leq f(u_k) = u_{k+1},
\]
since $f$ is non-decreasing. The result follows. \hfill$\blacksquare$
\end{pf}
\end{document}